\documentclass[a4paper,11pt]{amsart}

\usepackage{amssymb}

\usepackage[all]{xy}

\title[Homotopy branching space of flow]{The homotopy branching space of a flow}
\author[P. Gaucher]{Philippe Gaucher}
\address{Institut de Recherche Math\'ematique Avanc\'ee\\ ULP et
CNRS\\ 7 rue Ren\'e-Descartes\\ 67084 Strasbourg Cedex\\ France}
\email{gaucher@math.u-strasbg.fr}
\urladdr{http://www-irma.u-strasbg.fr/\~{}gaucher/}

\newcommand{\C}{\mathcal{C}}
\newcommand{\D}{\mathcal{D}}
\newcommand{\Z}{\mathbb{Z}}

\newcommand{\de}{\partial}
\newcommand{\p}\times
\renewcommand{\vec}{\overrightarrow}
\renewcommand{\P}{\mathbb{P}}

\newcommand{\be}{\begin{equation}}
\newcommand{\ee}{\end{equation}}
\newcommand{\bea}{\begin{eqnarray}}
\newcommand{\eea}{\end{eqnarray}}
\newcommand{\beas}{\begin{eqnarray*}}
\newcommand{\eeas}{\end{eqnarray*}}

\newtheorem{thm}{Theorem}[section]
\newtheorem{prop}[thm]{Proposition}

\newtheorem{rem}[thm]{Remark}

\newtheorem{defn}[thm]{Definition}
\newtheorem{propdef}[thm]{Proposition et Definition}

\newcommand{\bd}{\begin{defn}}
\newcommand{\ed}{\end{defn}}
\newcommand{\bcd}{\begin{defn}}
\newcommand{\ecd}{\end{defn}}
\newcommand{\bex}{\begin{exmp}}
\newcommand{\eex}{\end{exmp}}
\newcommand{\bp}{\begin{prop}}
\newcommand{\ep}{\end{prop}}
\newcommand{\bth}{\begin{thm}}
\renewcommand{\eth}{\end{thm}}
\newcommand{\br}{\begin{rem}}
\newcommand{\er}{\end{rem}}
\newcommand{\bpf}{\begin{proof}}
\newcommand{\epf}{\end{proof}}

\newcommand{\fl}[1]{\ar@{->}[l]_{#1}}
\newcommand{\fr}[1]{\ar@{->}[r]^-{#1}}
\newcommand{\fd}[1]{\ar@{->}[d]_{#1}}
\newcommand{\fu}[1]{\ar@{->}[u]^{#1}}
\newcommand{\f}[2]{\ar@{->}[#1]|{#2}}
\newcommand{\ff}[2]{\ar@2{->}[#1]|{#2}}
\newcommand{\frr}[1]{\ar@{->}[rr]^{#1}}

\renewcommand{\top}{{\mathbf{Top}}}

\newcommand{\ho}{{\mathbf{Ho}}}

\newcommand{\iso}{\cong}
\newcommand{\lp}{\left(}
\newcommand{\rp}{\right)}

\newcommand{\ot}{\otimes}
\newcommand{\eqc}[2]{\xymatrix{#1 \fr{\simeq}&#2}}

\newcommand{\vI}{\vec{I}}

\renewcommand{\geq}{\geqslant}

\def\cartesien{%
  \ar@{-}[]+R+<6pt,-2pt>;[]+RD+<6pt,-6pt>%
  \ar@{-}[]+D+<2pt,-6pt>;[]+RD+<6pt,-6pt>%
}
\def\cocartesien{%
  \ar@{-}[]+L+<-6pt,+2pt>;[]+LU+<-6pt,+6pt>%
  \ar@{-}[]+U+<-2pt,+6pt>;[]+LU+<-6pt,+6pt>%
}

\newcommand{\brm}[1]{\rm{\mathbf{#1}}}

\renewcommand{\top}{{\brm{Top}}}

\newcommand{\dtop}{{\brm{Flow}}}

\newcommand{\ttop}{{\brm{TOP}}}

\newcommand{\hda}{{I^{gl}_+\brm{cell}}}
\newcommand{\glob}{{\rm{Glob}}}

\DeclareMathOperator{\map}{Map}
\newcommand{\sis}{{\brm{SSet}}}

\newcommand{\liminj}{\varinjlim}

\makeatletter

\def\varholim@#1#2{%
  \vtop{\m@th\ialign{##\cr
    \hfil$#1\operator@font holim$\hfil\cr
    \noalign{\nointerlineskip\kern1.5\ex@}#2\cr
    \noalign{\nointerlineskip\kern-\ex@}\cr}}%
}
\def\holimproj{%
  \mathop{\mathpalette\varholim@{\leftarrowfill@\textstyle}}\nmlimits@
}
\def\holiminj{%
  \mathop{\mathpalette\varholim@{\rightarrowfill@\textstyle}}\nmlimits@
}

\newskip\@bigflushglue \@bigflushglue = -100pt plus 1fil

\def\bigcentering{\let\\\@centercr\rightskip\@bigflushglue%
\leftskip\@bigflushglue
\parindent\z@\parfillskip\z@skip}

\makeatother

\DeclareMathOperator{\hop}{ho\P}

\DeclareMathOperator{\id}{Id}

\begin{document}

\begin{abstract}

In this talk, I will explain the importance of the homotopy branching
space functor (and of the homotopy merging space functor) in
dihomotopy theory. 

\end{abstract}

\maketitle

\tableofcontents

\section{Introduction}

In \cite{survol}, the reader will be able to find a survey of the
different geometric approaches of concurrency. The model category of
flows was introduced in \cite{flow} to model higher dimensional
automata (HDA). It allows the study of HDA up to homotopy (cf. also
\cite{pgnote1,pgnote2}). A good notion of homotopy of flows must preserve 
the computer scientific properties of the HDA to be modeled like the
initial and final states, the deadlocks and the unreachable states.
In particular, it must preserve the direction of time, hence the
terminology
\textit{dihomotopy} for a contraction of \textit{directed homotopy}.  
This way, instead of working in the category of flows itself, one can
work in the localization of the category of flows with respect to
\textit{dihomotopy equivalences}.

I will explain in this talk the powerfulness of the homotopy branching
space functor in dihomotopy theory. The corresponding papers are
``Homotopy branching space and weak dihomotopy'' \cite{hobranching}
and ``A long exact sequence for the branching homology''
\cite{exbranching}.

\section{Model category}

If $\C$ is a category, one denotes by $Map(\C)$ the category whose
objects are the morphisms of $\C$ and whose morphisms are the commutative
squares of $\C$.

In a category $\C$, an object $x$ is \textit{a retract} of an object
$y$ if there exists $f:x\longrightarrow y$ and $g:y\longrightarrow x$
of $\C$ such that $g\circ f=\id_x$. A
\textit{functorial factorization} $(\alpha,\beta)$ of $\C$ is a pair of
functors from $Map(\C)$ to $Map(\C)$ such that for any $f$ object of
$Map(\C)$, $f=\beta(f)\circ \alpha(f)$.

\bd\cite{hovey,hir}
Let $i:A\longrightarrow B$ and $p:X\longrightarrow Y$ be maps in a
category $\C$. Then $i$ has the {\rm left lifting property} (LLP) with
respect to $p$ (or $p$ has the {\rm right lifting property} (RLP) with
respect to $i$) if for any commutative square
\[
\xymatrix{
A\fd{i} \fr{\alpha} & X \fd{p} \\
B \ar@{-->}[ru]^{g}\fr{\beta} & Y}
\]
there exists $g$ making both triangles commutative.  \ed

There are several versions of the notion of \textit{model\index{model
category} category}. The following definitions give the one we are
going to use.

\bd \cite{hovey,hir}
A {\rm model structure} on a category $\C$ is three subcategories of
$Map(\C)$ called weak equivalences, cofibrations, and fibrations, and
two functorial factorizations $(\alpha,\beta)$ and $(\gamma,\delta)$
satisfying the following properties\thinspace:
\begin{enumerate}
\item (2-out-of-3) If $f$ and $g$ are morphisms of $\C$ such that $g\circ f$
is defined and two of $f$, $g$ and $g\circ f$ are weak equivalences,
then so is the third.
\item (Retracts)
If $f$ and $g$ are morphisms of $\C$ such that $f$ is a retract
of $g$ and $g$ is a weak equivalence, cofibration, or fibration, then
so is $f$.
\item (Lifting)
Define a map to be a trivial cofibration if it is both a cofibration
and a weak equivalence. Similarly, define a map to be a trivial
fibration if it is both a fibration and a weak equivalence. Then
trivial cofibrations have the LLP\index{left lifting property} with
respect to fibrations, and cofibrations have the LLP\index{left
lifting property} with respect to trivial fibrations.
\item (Factorization) For any morphism $f$, $\alpha(f)$ is a cofibration,
$\beta(f)$ a trivial fibration, $\gamma(f)$ is a trivial cofibration ,
and $\delta(f)$ is a fibration.
\end{enumerate}
\ed

\bd\cite{hovey,hir}
A {\rm model\index{model category} category} is a complete and
cocomplete category $\C$ together with a model structure on $\C$. \ed

\begin{propdef}\cite{hovey,hir}
A Quillen adjunction is a pair of adjoint functors
$F:\C\rightleftarrows \D:G$ between the model categories $\C$ and $\D$
such that one of the following equivalent properties holds\thinspace:
\begin{enumerate}
\item if $f$ is a cofibration (resp. a trivial cofibration), then so
does $F(f)$
\item if $g$ is a fibration (resp. a trivial fibration), then so
does $G(g)$.
\end{enumerate}
One says that $F$ is a {\rm left Quillen functor}.
One says that $G$ is a {\rm right Quillen functor}.
\end{propdef}

\bd\cite{hovey,hir}
An object $X$ of a model category $\C$ is {\rm cofibrant}
(resp. {\rm fibrant}) if and only if the canonical morphism
$\varnothing\longrightarrow X$ from the initial object of $\C$ to $X$
(resp. the canonical morphism $X\longrightarrow \mathbf{1}$ from $X$
to the final object $\mathbf{1}$) is a cofibration (resp.  a
fibration). \ed

For any object $X$ of a model category, the canonical morphism 
$\varnothing_X:\varnothing\longrightarrow X$ from the initial object to $X$ can be 
factored as a composite 
\[\xymatrix{\varnothing\ar@{->}[rr]^-{\alpha(\varnothing_X)}&& Q(X)\ar@{->}[rr]^-{\beta(\varnothing_X)}&& X}\]
where, by definition, $Q(X)$ is a cofibrant object which is weakly
equivalent to $X$. The functor $Q:\C\longrightarrow \C$ is called the
\textit{cofibrant replacement functor}.

\section{Reminder about the category of flows}

In the sequel, any topological space will be supposed to be compactly
generated (more details for this kind of topological spaces in
\cite{MR90k:54001,MR2000h:55002}, the appendix of \cite{Ref_wH} and
also the preliminaries of
\cite{flow}).

Let $n\geq 1$. Let $\mathbf{D}^n$ be the closed $n$-dimensional disk.
Let $\mathbf{S}^{n-1}=\de \mathbf{D}^n$ be the boundary of
$\mathbf{D}^n$ for $n\geq 1$. Notice that $\mathbf{S}^0$ is the
discrete\index{discrete topology} two-point topological space
$\{-1,+1\}$.  Let $\mathbf{D}^0$ be the one-point topological space.
Let $\mathbf{S}^{-1}=\varnothing$ be the empty set. The following 
theorem is well-known.

\bth \label{usualmodel}
\cite{hir,hovey} The category of compactly generated topological spaces
$\top$ can be given a model structure such that:
\begin{enumerate}
\item The weak equivalences are the weak homotopy equivalences.
\item The fibrations (sometime called Serre fibrations)
are the continuous maps satisfying the RLP (right lifting property) with respect
to the continuous maps $\mathbf{D}^n\longrightarrow [0,1]\p \mathbf{D}^n$ such that $x\mapsto (0,x)$ and
for $n\geq 0$.
\item The cofibrations are the continuous maps satisfying the LLP (left lifting property) with
respect to any maps satisfying the RLP with respect to the
inclusion  maps $\mathbf{S}^{n-1}\longrightarrow \mathbf{D}^n$.
\item Any topological space is fibrant.
\item The homotopy equivalences arising from this model structure coincide
with the usual one.
\end{enumerate}
\eth

\bd \cite{flow}
A {\rm flow} $X$ consists of a topological space $\P X$, a discrete
space $X^0$, two continuous maps $s$ and $t$ from $\P X$ to $X^0$ and
a continuous and associative map $*:\{(x,y)\in \P X\p \P X;
t(x)=s(y)\}\longrightarrow \P X$ such that $s(x*y)=s(x)$ and
$t(x*y)=t(y)$.  A morphism of flows $f:X\longrightarrow Y$ consists of
a set map $f^0:X^0\longrightarrow Y^0$ together with a continuous map
$\P f:\P X\longrightarrow \P Y$ such that $f(s(x))=s(f(x))$,
$f(t(x))=t(f(x))$ and $f(x*y)=f(x)*f(y)$. The corresponding category
will be denoted by $\dtop$.
\ed

The topological space $X^0$ is called the \textit{$0$-skeleton} of
$X$. The topological space $\P X$ is called the \textit{path space}
and its elements the \textit{non constant execution paths} of $X$. The
initial object $\varnothing$ of $\dtop$ is the empty set. The terminal
object $\mathbf{1}$ is the flow defined by $\mathbf{1}^0=\{0\}$, $\P
\mathbf{1}=\{u\}$ and necessarily $u*u=u$.

\bd\cite{flow}
Let $Z$ be a topological space. Then the {\rm globe} of $Z$ is the
flow $\glob(Z)$ defined as follows: $\glob(Z)^0=\{0,1\}$,
$\P\glob(Z)=Z$, $s=0$, $t=1$ and the composition law is trivial. \ed

\bth \cite{flow}\label{model}  The category of flows
can be given a model structure such that:
\begin{enumerate}
\item The weak equivalences  are the {\rm weak S-homotopy equivalences}, that is
a morphism of flows $f:X\longrightarrow Y$ such that $f:X^0\longrightarrow Y^0$ is an
isomorphism of sets and $f:\P X\longrightarrow \P Y$ a weak homotopy
equivalence of topological spaces.
\item The fibrations are the continuous maps satisfying the RLP
 with respect
to the morphisms $\glob(\mathbf{D}^n)\longrightarrow \glob([0,1]\p \mathbf{D}^n)$ for $n\geq
0$. The fibrations are exactly the morphisms of flows $f:X\longrightarrow
Y$ such that $\P f:\P X\longrightarrow \P Y$ is a Serre fibration of $\top$.
\item The cofibrations are the morphisms satisfying the LLP  with
respect to any map satisfying the RLP with respect to the morphisms
$\glob(\mathbf{S}^{n-1})\longrightarrow \glob(\mathbf{D}^n)$  for
$n\geq 0$ and with respect to the morphisms $\varnothing\longrightarrow
\{0\}$ and $\{0,1\}\longrightarrow \{0\}$.
\item Any flow is fibrant.
\end{enumerate}
\eth

Let $I^{gl}$ be the set of morphisms of flows
$\glob(\mathbf{S}^{n-1})\rightarrow \glob(\mathbf{D}^n)$ for $n\geq 0$. Denote by
$I^{gl}_+$ be the union of $I^{gl}$ with the two morphisms of flows
$R:\{0,1\}\rightarrow \{0\}$ and $C:\varnothing\subset \{0\}$.

\bd\label{cell}\cite{hobranching} 
An {\rm $I^{gl}_+$-cell complex} is a flow $X$ such that the canonical
morphism of flows $\varnothing\longrightarrow X$ from the initial object
of $\dtop$ to $X$ is a transfinite composition of pushouts of elements
of $I^{gl}_+$. The full and faithful subcategory of $\dtop$ whose
objects are the $I^{gl}_+$-cell complexes will be denoted by $\hda$.
\ed

The category $\hda$ of $I^{gl}_+$-cell complexes is a subcategory of
the category of flows which is sufficient to model higher dimensional
automata (HDA), at least those modeled by precubical sets
\cite{diCW,cridlig96implementing}. This geometric model of HDA is
designed to define and study equivalence relations preserving the
compu\-ter-scientific properties of the HDA to be modeled so that it
then suffices to work in convenient localizations of $\hda$. The
properties which are preserved are for instance the initial or final
states, the presence or not of deadlocks and of unreachable states
\cite{flow}.

The cofibrant replacement functor is a functor $Q:\dtop\longrightarrow
\hda$.  The flows coming from concrete HDAs are all cofibrant.

\section{The homotopy branching space functor}

The branching space of a flow is the space of germs of non-constant
execution paths beginning in the same way. The branching space functor
$\P^-$ from the category of flows $\dtop$ to the category of compactly
generated topological spaces $\top$ was also introduced in \cite{flow}
to fit the definition of the branching semi-globular nerve of a strict
globular $\omega$-category modeling an HDA introduced in
\cite{fibrantcoin}.

\bp\cite{flow,hobranching} \label{universalpm} 
Let $X$ be a flow. There exists a topological space
$\P^-X$ unique up to homeomorphism and a continuous map $h^-:\P
X\longrightarrow \P^- X$ satisfying the following universal
property:
\begin{enumerate}
\item For any $x$ and $y$ in $\P X$ such that $t(x)=s(y)$, the equality
$h^-(x)=h^-(x*y)$ holds.
\item Let $\phi:\P X\longrightarrow Y$ be a
continuous map such that for any $x$ and $y$ of $\P X$ such that
$t(x)=s(y)$, the equality $\phi(x)=\phi(x*y)$ holds. Then there exists a
unique continuous map $\overline{\phi}:\P^-X\longrightarrow Y$ such that
$\phi=\overline{\phi}\circ h^-$.
\end{enumerate}
Moreover, one has the homeomorphism
\[\P^-X\iso \bigsqcup_{\alpha\in X^0} \P^-_\alpha X\]
where $\P^-_\alpha X:=h^-\lp \bigsqcup_{\beta\in
X^0} \P_{\alpha,\beta} X\rp$. The mapping $X\mapsto \P^-X$
yields a functor $\P^-$ from $\dtop$ to $\top$. \ep

\bd\cite{flow,hobranching} 
Let $X$ be a flow. The topological space $\P^-X$ is called the {\rm
branching space} of the flow $X$. \ed

\bp \cite{hobranching} \label{contre1}
There exists a weak S-homotopy equivalence of flows
$f:X\longrightarrow Y$ such that the topological spaces $\P^-X$ and
$\P^-Y$ are not weakly homotopy equivalent.
\ep

The idea for the proof of Proposition~\ref{contre1} is as follows. For
a given flow $X$, by Proposition~\ref{universalpm}, the topological space 
$\P^- X$ is the coequalizer of the continuous map $\P X\p_{X^0} \P X\longrightarrow 
\P X$ induced by the composition law of $X$ and of the projection map $\P X\p_{X^0} \P X\longrightarrow 
\P X$ on the first factor. And one cannot expect a coequalizer to transform 
a objectwise weak homotopy equivalence into a weak homotopy
equivalence. One must use a kind of homotopy coequalizer instead.

If two flows are weakly S-homotopy equivalent, then they are supposed
to satisfy the same computer-scientific properties. With the example
above, one obtains two such flows but with very different branching
spaces.  But

\bth \cite{hobranching} 
If $f:X\longrightarrow Y$ is a weak S-homotopy equivalence of flows
between cofibrant flows, then the topological spaces $\P^-X$ and
$\P^-Y$ are homotopy equivalent. \eth

This suggests that the definition of the branching space is the good
one up to homotopy for cofibrant flows. Indeed, we have the theorems:

\bth\cite{hobranching} \label{Qad}
There exists a functor $C^-:\top\longrightarrow \dtop$ such that the
pair of functors $\P^-:\dtop\rightleftarrows\top:C^-$ is a Quillen
adjunction. In particular, there is an homeomorphism $\P^-(\liminj X_i)\iso
\liminj \P^-X_i$. \eth

\bd
The {\rm homotopy branching space} $\hop^- X$ of a flow $X$ is by
definition the topological space $\P^-Q(X)$. \ed

\bth\cite{hobranching}\label{leftderived}
The functor $\hop^-:\dtop\longrightarrow \top\longrightarrow
\ho(\top)$ satisfies the following universal property: if
$F:\dtop\longrightarrow \ho(\top)$ is another functor sending weak
S-homotopy equivalences to isomorphisms and if there exists a natural
transformation $F\Rightarrow \P^-$, then the latter natural
transformation factors uniquely as a composite $F\Rightarrow \hop^-
\Rightarrow \P^-$. \eth

Up to homotopy, the homotopy bran\-ching space $\hop^-(X)$ is
well-defined and coincides with $\P^-X$ for any cofibrant flow, so in
particular for any flow coming from a HDA.  The behavior of the
branching space functor and the homotopy branching space functor are
the same up to homotopy for flows modeling HDAs and may differ for
other flows.

\section{The homotopy merging space functor}

This is the dual version of the preceding functor. Some results are
collected in this section about it.

\bp\cite{hobranching} \label{universalpp} Let $X$ be a flow. There exists a topological space
$\P^+X$ unique up to homeomorphism and a continuous map $h^+:\P
X\longrightarrow \P^+ X$ satisfying the following universal
property\thinspace:
\begin{enumerate}
\item For any $x$ and $y$ in $\P X$ such that $t(x)=s(y)$, the equality
$h^+(y)=h^+(x*y)$ holds.
\item Let $\phi:\P X\longrightarrow Y$ be a
continuous map such that for any $x$ and $y$ of $\P X$ such that
$t(x)=s(y)$, the equality $\phi(y)=\phi(x*y)$ holds. Then there exists a
unique continuous map $\overline{\phi}:\P^+X\longrightarrow Y$ such that
$\phi=\overline{\phi}\circ h^+$.
\end{enumerate}
Moreover, one has the homeomorphism
\[\P^+X\iso \bigsqcup_{\alpha\in X^0} \P^+_\alpha X\]
where $\P^+_\alpha X:=h^+\lp \bigsqcup_{\beta\in X^0}
\P_{\beta,\alpha} X\rp$. The mapping $X\mapsto \P^+X$ yields a
functor $\P^+:\dtop\longrightarrow \top$. \ep

\bd\cite{hobranching} Let $X$ be a flow. The topological space $\P^+X$ is called the {\rm merging
space} of the flow $X$. \ed

\bth\cite{hobranching} \label{Qadp}
There exists a functor $C^+:\top\longrightarrow \dtop$
such that the pair of functors $\P^+:\dtop\rightleftarrows
\top:C^+$ is a Quillen adjunction. In particular, there is an  homeomorphism
$\P^+(\liminj X_i)\iso \liminj \P^+X_i$. \eth

\bd\cite{hobranching} The {\rm homotopy merging space} $\hop^+ X$ of a flow $X$ is by definition
the topological space $\P^+Q(X)$. \ed

\bth\cite{hobranching}
The functor $\hop^+:\dtop\longrightarrow \top\longrightarrow
\ho(\top)$ satisfies the following universal property\thinspace: if
$F:\dtop\longrightarrow \ho(\top)$ is another functor sending weak
S-homotopy equivalences to isomorphisms and if there exists a natural
transformation $F\Rightarrow \P^+$, then the latter natural
transformation factors uniquely as a composite $F\Rightarrow \hop^+
\Rightarrow \P^+$. \eth

\section{First application: studying weak dihomotopy}

The class $\mathcal{S}$ of weak S-homotopy equivalences is an example
of class of morphisms of flows which is supposed to preserve various
computer-scientific properties. This class of morphisms of flows
satisfies the following properties:
\begin{enumerate}
\item The two-out-of-three axiom, that is if two of the three morphisms
$f$, $g$ and $g\circ f$ belong to $\mathcal{S}$, then so does the
third one: this condition means that the class $\mathcal{S}$
defines an equivalence relation.
\item The embedding functor $I:\hda\longrightarrow \dtop$ induces
a functor $\overline{I}:\hda[\mathcal{S}^{-1}]\longrightarrow
\dtop[\mathcal{S}^{-1}]$ between the localization of respectively the category
of $I^{gl}_+$-cell complexes and the category of flows with respect to
weak S-homotopy equivalences which is an equivalence of categories. In
particular, it reflects isomorphisms, that is $X\iso Y$ if and only if
$\overline{I}(X)\iso \overline{I}(Y)$. In this case, one can use the
whole category of flows which is a richer mathematical framework.
\end{enumerate}

The class of T-homotopy equivalences was introduced in \cite{flow} to
identify $I^{gl}_+$-cell complexes equivalent from a
computer-scientific viewpoint and which are not identified in
$\hda[\mathcal{S}^{-1}]$. Indeed, if two objects $X$ and $Y$ of
$\hda[\mathcal{S}^{-1}]$ are isomorphic, then the $0$-skeletons $X^0$
and $Y^0$ are isomorphic.  The merging of the notions of weak
S-homotopy equivalence and T-homotopy equivalence yields the class
$\mathcal{ST}_0$ of $0$-dihomotopy equivalences.

\bd \cite{flow} Let $X$ be a flow. Let $A$ and $B$ be two subsets of $X^0$. One says that
$A$ is {\rm surrounded} by $B$ (in $X$) if for any $\alpha\in A$,
either $\alpha \in B$ or there exists execution paths $\gamma_1$ and
$\gamma_2$ of $\P X$ such that $s(\gamma_1)\in B$,
$t(\gamma_1)=s(\gamma_2)=\alpha$ and $t(\gamma_2)\in B$. We denote
this situation by $A\lll B$.  \ed

\bd  \cite{flow} Let $X$ be a flow. Let $A$ be a subset of $X^0$. Then
the {\rm restriction} $X\!\restriction_A$ of $X$ over $A$ is the
unique flow such that $\left(X\!\restriction_A\right)^0=A$ and
\[\P \left(X\!\restriction_A\right) = \bigsqcup_{\left(\alpha,\beta\right)\in A\p A}
\P_{\alpha,\beta} X\] equipped with the topology induced by
the one of $\P X$. \ed

\bd\cite{flow}\label{def0}
A morphism of flows $f:X\longrightarrow Y$ is a {\rm $0$-dihomotopy
equivalence} if and only if the following conditions are
satisfied\thinspace:
\begin{enumerate}
\item The morphism of flows $f:X\longrightarrow
Y\!\restriction_{f(X^0)}$ is a weak S-homotopy equivalence of
flows. In particular, the set map $f^0:X^0\longrightarrow Y^0$ is
one-to-one.
\item For $\alpha\in Y^0\backslash f(X^0)$, the topological spaces
$\P^-_\alpha Y$ and $\P^+_\alpha Y$ are singletons.
\item $Y^0\lll f(X^0)$.
\end{enumerate}
The class of $0$-dihomotopy equivalences is denoted by $\mathcal{ST}_0$.
\ed

But it turns out that

\bth \cite{hobranching}
The functor $\hda[\mathcal{ST}_0^{-1}]\longrightarrow
\dtop[\mathcal{ST}_0^{-1}]$ does not reflect isomorphisms.
More precisely, there exists an $I^{gl}_+$-cell complex $\vec{C}_3$
corresponding to the concurrent execution of three calculations which
is not isomorphic in $\hda[\mathcal{ST}_0^{-1}]$ to the directed
segment $\vI$, although the same flow $\vec{C}_3$ is isomorphic to
$\vI$ in $\dtop[\mathcal{ST}_0^{-1}]$. 
\eth

The correct behavior is the one of $\mathcal{ST}_0$ in
$\dtop[\mathcal{ST}_0^{-1}]$.  Indeed, an HDA representing the
concurrent execution of $n$ processes must be equivalent to the
directed segment in a good homotopical approach of concurrency. The
interpretation of this fact is therefore that the class
$\mathcal{ST}_0$ of $0$-dihomotopy equivalences is not big enough.

\bd\label{def1}\cite{hobranching}
A morphism of flows $f:X\longrightarrow Y$ is a {\rm $1$-dihomotopy
equivalence} if and only if the following conditions are
satisfied\thinspace:
\begin{enumerate}
\item The morphism of flows $f:X\longrightarrow
Y\!\restriction_{f(X^0)}$ is a weak S-homotopy equivalence  of flows. In
particular, the set map $f^0:X^0\longrightarrow Y^0$ is one-to-one.
\item For $\alpha\in Y^0\backslash f(X^0)$, the topological spaces
$\P^-_\alpha Y$ and $\P^+_\alpha Y$ are weakly contractible.
\item $Y^0\lll f(X^0)$.
\end{enumerate}
The class of $1$-dihomotopy equivalences is denoted by $\mathcal{ST}_1$.
\ed

Any $0$-dihomotopy equivalence is of course a $1$-dihomotopy
equivalence. Moreover, the composite of a weak S-homotopy equivalence
with a T-homotopy equivalence can already give an element of
$\mathcal{ST}_1\backslash \mathcal{ST}_0$ ! And

\bth \cite{hobranching} 
By slightly weakening the notion of T-homotopy as above, one obtains a
class of morphisms $\mathcal{ST}_1$ with $\mathcal{ST}_0\subset
\mathcal{ST}_1$ and such that the flows $\vec{C}_3$ and $\vI$ become
isomorphic in the localization $\hda[\mathcal{ST}_1^{-1}]$. \eth

There are actually two natural ways of weakening the definition of
$\mathcal{ST}_0$. One can replace in the statement the word
\textit{singleton} either by the word \textit{weakly contractible}, or
by the word \textit{contractible}. This way, one obtains another class
of morphisms $\mathcal{ST'}_1$ with $\mathcal{ST'}_1\subset
\mathcal{ST}_1$ and one has:

\bth \cite{hobranching} 
The localizations $\hda[\mathcal{ST'}_1^{-1}]$ and
$\hda[\mathcal{ST}_1^{-1}]$ are equivalent.
\eth

Unfortunately, one has

\bp \cite{hobranching} 
The composite of two morphisms of $\mathcal{ST}_1$ does not
necessarily belong to $\mathcal{ST}_1$. \ep

Using the homotopy branching space functor, a new class
$\mathcal{ST}_2$ of morphisms of flows is introduced.

\bd\label{def2}\cite{hobranching}
A morphism of flows $f:X\longrightarrow Y$ is a {\rm $2$-dihomotopy
equivalence} if and only if the following conditions are
satisfied\thinspace:
\begin{enumerate}
\item The morphism of flows $f:X\longrightarrow
Y\!\restriction_{f(X^0)}$ is a weak S-homotopy equivalence  of flows. In
particular, the set map $f^0:X^0\longrightarrow Y^0$ is one-to-one.
\item For $\alpha\in Y^0\backslash f(X^0)$, the topological spaces
$\hop^-_\alpha Y$ and $\hop^+_\alpha Y$ are weakly contractible.
\item $Y^0\lll f(X^0)$.
\end{enumerate}
The class of $2$-dihomotopy equivalences is denoted by $\mathcal{ST}_2$.
\ed

And:

\bth  \cite{hobranching}
One has the equivalence of categories
\[\eqc{\hda[\mathcal{ST}_1^{-1}]}{\hda[\mathcal{ST}_2^{-1}]}\]
where $\hda[\mathcal{ST}_1^{-1}]$ (resp. $\hda[\mathcal{ST}_2^{-1}]$)
is the localization of the category of $I^{gl}_+$-cell complexes with
respect to $1$-dihomotopy equivalences (resp.  $2$-dihomotopy
equivalences). $\mathcal{ST}_2$ is closed under composition.  Moreover
the embedding functor $I:\hda\longrightarrow \dtop$ induces an
equivalence of categories
\[\overline{I}:\eqc{\hda[\mathcal{ST}_2^{-1}]}{\dtop[\mathcal{ST}_2^{-1}]}.\]
In particular, the functor $\hda[\mathcal{ST}_2^{-1}]\longrightarrow
\dtop[\mathcal{ST}_2^{-1}]$ reflects isomorphisms. 
\eth

The property $f\in\mathcal{ST}_2\hbox{ and }g\circ f\in
\mathcal{ST}_2\Longrightarrow g\in \mathcal{ST}_2$ has no reasons to
be satisfied by $2$-dihomotopy equivalences.  Indeed, if both $g\circ
f$ and $f$ are two one-to-one set maps, then $g$ has no reasons to be
one-to-one as well. Therefore in order to understand the isomorphisms
of $\dtop[\mathcal{ST}_2^{-1}]$, we may introduce another
construction.

\bd\cite{hobranching}
Let $X$ be a flow. Then a subset $A$ of $X^0$ is {\rm essential} if
$X^0\lll A$ and if for any $\alpha\notin A$, both topological spaces
$\hop^-_\alpha X$ and $\hop^+_\alpha X$ are weakly contractible. \ed

\bd\cite{hobranching}
A morphism of flows $f:X\longrightarrow Y$ is a {\rm $3$-dihomotopy
equivalence} if the following conditions are satisfied\thinspace:
\begin{enumerate}
\item $A\subset X^0$ is essential if and only if $f(A)\subset Y^0$ is
essential
\item for any essential $A\subset X^0$ there exists an essential subset
$B\subset A$ such that the restriction $f:X\!\restriction_B\longrightarrow Y\!\restriction_{f(B)}$ is a weak S-homotopy equivalence.
\end{enumerate}
The class of $3$-dihomotopy equivalences is denoted by $\mathcal{ST}_3$.
\ed

\bth \cite{hobranching}
The localizations $\hda[\mathcal{ST}_2^{-1}]$ and
$\hda[\mathcal{ST}_3^{-1}]$ are equivalent and the class of morphisms
$\mathcal{ST}_3$ satisfies the two-out-of-three axiom.  Moreover the
embedding functor $I:\hda\longrightarrow \dtop$ induces an equivalence
of categories
\[\overline{I}:\eqc{\hda[\mathcal{ST}_3^{-1}]}{\dtop[\mathcal{ST}_3^{-1}]}.\]
In particular, the functor $\hda[\mathcal{ST}_3^{-1}]\longrightarrow
\dtop[\mathcal{ST}_3^{-1}]$ reflects isomorphisms. 
\eth

The class $\mathcal{ST}_2$ does not satisfy the two-out-of-three axiom
but is invariant by retract. The class $\mathcal{ST}_3$ does satisfy
the two-out-of-three axiom but is probably not invariant by
retract. So none of the definitions above allows to describe the
isomorphisms of $\hda[\mathcal{ST}_2^{-1}]$. The situation can be summarized 
with the following diagram: 

\begin{center}
{
\xymatrix{
&& \hda \ar@/_10pt/[lld]\ar@{->}[ld]\ar@/^5pt/[rd]\fd{}&&\\
\hda[\mathcal{S}^{-1}]\fd{\simeq}\ar@{->}[r]^{\not\simeq}&\fd{\not\simeq}\hda[\mathcal{ST}_0^{-1}]
\ar@{->}[r]^{\not\simeq}&\hda[\mathcal{ST}_1^{-1}]\fd{\simeq ??}\fr{\simeq}& \hda[\mathcal{ST}_2^{-1}]\simeq
\hda[\mathcal{ST}_3^{-1}]\fd{\simeq}\\
\dtop[\mathcal{S}^{-1}]\ar@{->}[r]^{\not\simeq}&\dtop[\mathcal{ST}_0^{-1}]
\ar@{->}[r]^{\simeq??}&\dtop[\mathcal{ST}_1^{-1}]&\ar@{->}[l]_-{\simeq ??}\dtop[\mathcal{ST}_2^{-1}]\simeq
\dtop[\mathcal{ST}_3^{-1}]\\
&& \dtop \ar@/^10pt/[llu]\ar@{->}[lu]\ar@{->}[u]\ar@/_5pt/[ru]&&
}
}\end{center}

The symbol $\simeq??$ means that we do not know whether the functor is
an equivalence of categories or not. The symbol $\not\simeq$ means
that the corresponding functor is not an equivalence.

\section{Second application: a long exact sequence for the branching homology}

The category of flows is a simplicial model category
\cite{exbranching} in the following sense:

\bd \cite{Quillen,hovey,hir} \label{cidessus}x
A {\rm simplicial model category} is a model category $\C$ together with a simplicial
set $\map(X,Y)$ for any object $X$ and $Y$ of $\C$ satisfying the following
axioms:
\begin{enumerate}
\item the set $\map(X,Y)_0$ is canonically isomorphic to $\C(X,Y)$
\item for any object $X$, $Y$ and $Z$, there is a morphism of simplicial
sets \[\map(Y,Z)\p \map(X,Y)\longrightarrow \map(X,Z)\] which is associative
\item for any object $X$ of $\C$ and any simplicial set $K$, there exists
an object $X\ot K$  of $\C$ such
that there exists a natural isomorphism of simplicial sets \[\map(X\ot
K,Y)\iso \map (K,\map(X,Y))\]
\item for any object $X$ of $\C$ and any simplicial set $K$, there exists
an object $X^K$ such
that there exists a natural isomorphism of simplicial sets
\[\map(X,Y^K)\iso \map (K,\map(X,Y))\]
\item for any cofibration $i:A\longrightarrow B$ and any fibration $p:X\longrightarrow Y$
of $\C$, the morphism of simplicial sets
\[Q(i,p):\map(B,X)\longrightarrow \map(A,X)\p_{\map(A,Y)}\map(B,Y)\]
is a fibration of simplicial sets. Moreover if either $i$ or $p$ is trivial,
then the fibration $Q(i,p)$ is trivial as well.
\end{enumerate}
\ed

Recall that there exists a
pair of adjoint functors $|-|:\sis \rightleftarrows \top : S_*$ where
$|-|$ is the geometric realization functor and $S_*$ the singular
nerve functor. The $n$-simplex of $\sis$ is denoted by $\Delta[n]$.
Its boundary is denoted by $\de\Delta[n-1]$.
Let $\Delta^n$ be the $n$-dimensional simplex.

The category of compactly generated topological spaces $\top$ is a
simplicial model category by setting $\map(X,Y)_n:=\top(X\p
\Delta^n,Y)$, $X\ot K:= X\p |K|$ and $X^K:=\ttop(|K|,X)$. The category
of simplicial sets $\sis$ is a simplicial model category as well by
setting $\map(X,Y)_n:=\top(X\p \Delta[n],Y)$, $X\ot K:=X\p K$ and
$X^K:=\map(K,X)$ \cite{Quillen}.

This means that the model category of flows can be enriched over the 
category of simplicial sets and that the enrichment is compatible with the model structure in the sense 
of Definition~\ref{cidessus}. The symbol $\Delta^n$ is the simplicial set 
corresponding to the $n$-dimensional simplex.

Because of the existence of this enrichment, there exist explicit
formulae for homotopy colimits \cite{hir}.  In particular, the
homotopy pushout of a diagram of flows looks as follows:

\bd  \cite{hir} The homotopy pushout of the diagram of flows 
\[
\xymatrix{
A\fr{}\fd{}  & B\\
C  & }
\]
is the colimit of the diagram of flows 
\[
\xymatrix{
 & A\ot \Delta^0\fd{} \fr{} & B \\
A\ot \Delta^0\fd{} \fr{} & A\ot \Delta^1 & \\
C && }
\]
\ed

It is then very easy to prove the:

\bth\cite{exbranching} \label{preho}
Let $X$ be a diagram of flows. Then the topological spaces $\holiminj
\hop^-(X)$ and $\hop^- (\holiminj X)$ are homotopy equivalent (they
are both cofibrant indeed). So in particular, the homotopy branching space 
functor commutes with homotopy pushouts. 
\eth

\bd\cite{exbranching} Let $f:X\longrightarrow Y$ be a morphism of flows. The {\rm cone} $Cf$ of $f$
is the homotopy pushout in the category of flows
\[
\xymatrix{
X\fr{f} \fd{} & Y \fd{} \\
\textbf{1} \fr{} & Cf
}
\]
where $\textbf{1}$ is the terminal flow.
\ed

From the theorem

\bth\cite{exbranching}\label{point} The homotopy branching space of the terminal
flow is contractible. \eth

one can easily deduce a long exact sequence for the branching
homology. 

\bd\cite{exbranching} Let $X$ be a flow. Then the $(n+1)$st branching homology group
$H_{n+1}^-(X)$ is defined as the  $n$st homology group of the augmented
simplicial set $\mathcal{N}^-_*(X)$ defined as follows:
\begin{enumerate}
\item $\mathcal{N}^-_n(X)=S_n(\hop^-X)$ for $n\geq 0$
\item $\mathcal{N}^-_{-1}(X)=X^0$
\item the augmentation map $\epsilon:S_0(\hop^-X)\longrightarrow X^0$
is induced by the mapping $\gamma\mapsto s(\gamma)$ from $\hop^-X=S_0(\hop^-X)$
to $X^0$.
\end{enumerate}
\ed

\bth\cite{exbranching} For any flow $X$, one has
\begin{enumerate}
\item $H_0^-(X)=\Z X^0/Im(s)$
\item the short exact sequence $0\rightarrow H_1^-(X)\rightarrow H_0(\hop^-X)\rightarrow \Z \hop^-X/Ker(s)\rightarrow 0$
\item $H_{n+1}^-(X)=H_n(\hop^-X)$ for $n\geq 1$.
\end{enumerate}
\eth

\bth\cite{exbranching} 
For any morphism of flows $f:X\longrightarrow Y$, one has the long
exact sequence
\beas
&& \dots \rightarrow H_{n}^-(X) \rightarrow H_{n}^-(Y) \rightarrow H_{n}^-(Cf)\rightarrow  \dots \\
&& \dots \rightarrow H_{3}^-(X) \rightarrow H_{3}^-(Y) \rightarrow H_{3}^-(Cf)\rightarrow \\
&& H_{2}^-(X) \rightarrow H_{2}^-(Y) \rightarrow H_{2}^-(Cf)\rightarrow \\
&& H_0(\hop^-X) \rightarrow H_0(\hop^-Y)\rightarrow  H_0(\hop^- Cf)\rightarrow 0.
\eeas
\eth

The functors $X\mapsto H_{n}^-(X)$ for $n\geq 0$ are invariant up to $2$-dihomotopy
equivalence. The functor $X\mapsto H_0(\hop^-X)$ is only invariant up
to weak S-homotopy equivalence. So the long exact sequence above is
not satisfactory. It still remains to find an exact sequence whose
each term would be a functor invariant up to $2$-dihomotopy
equivalence.


\begin{thebibliography}{10}

\bibitem{MR90k:54001}
R.~Brown.
\newblock {\em Topology}.
\newblock Ellis Horwood Ltd., Chichester, second edition, 1988.
\newblock A geometric account of general topology, homotopy types and the
  fundamental groupoid.

\bibitem{cridlig96implementing}
R.~Cridlig.
\newblock Implementing a static analyzer of concurrent programs: Problems and
  perspectives.
\newblock In {\em Logical and Operational Methods in the Analysis of Programs
  and Systems}, pages 244--259, 1996.

\bibitem{flow}
P.~Gaucher.
\newblock {A Convenient Category for The Homotopy Theory of Concurrency}.
\newblock arXiv:math.AT/0201252.

\bibitem{exbranching}
P.~Gaucher.
\newblock {A long exact sequence for the branching homology}.
\newblock arXiv:math.AT/0305169.

\bibitem{hobranching}
P.~Gaucher.
\newblock {Homotopy branching space and weak dihomotopy}.
\newblock arXiv:math.AT/0304112.

\bibitem{fibrantcoin}
P.~Gaucher.
\newblock The branching nerve of {HDA} and the {K}an condition.
\newblock {\em Theory and Applications of Categories}, 11(3):pp.75--106, 2003.

\bibitem{pgnote1}
P.~Gaucher.
\newblock Concurrent process up to homotopy ({I}).
\newblock {\em C. R. Acad. Sci. Paris Ser. I Math.}, 336(7):593--596, 2003.
\newblock French.

\bibitem{pgnote2}
P.~Gaucher.
\newblock Concurrent process up to homotopy ({II}).
\newblock {\em C. R. Acad. Sci. Paris Ser. I Math.}, 336(8):647--650, 2003.
\newblock French.

\bibitem{diCW}
P.~Gaucher and E.~Goubault.
\newblock Topological deformation of higher dimensional automata.
\newblock {\em Homology, Homotopy and Applications}, 5(2):p.39--82, 2003.

\bibitem{survol}
E.~Goubault.
\newblock Some geometric perspectives in concurrency theory.
\newblock {\em Homology, Homotopy and Applications}, 5(2):p.95--136, 2003.

\bibitem{hir}
P.~S. Hirschhorn.
\newblock {\em Model categories and their localizations}, volume~99 of {\em
  Mathematical Surveys and Monographs}.
\newblock American Mathematical Society, Providence, RI, 2003.

\bibitem{hovey}
M.~Hovey.
\newblock {\em Model categories}.
\newblock American Mathematical Society, Providence, RI, 1999.

\bibitem{Ref_wH}
L.~G. Lewis.
\newblock {\em The stable category and generalized Thom spectra}.
\newblock PhD thesis, University of Chicago, 1978.

\bibitem{MR2000h:55002}
J.~P. May.
\newblock {\em A concise course in algebraic topology}.
\newblock University of Chicago Press, Chicago, IL, 1999.

\bibitem{Quillen}
D.~G. Quillen.
\newblock {\em Homotopical algebra}.
\newblock Lecture Notes in Mathematics, No. 43. Springer-Verlag, Berlin, 1967.

\end{thebibliography}
\end{document}